\magnification=\magstephalf
\input amstex
\loadbold
\documentstyle{amsppt}
\refstyle{A}
\NoBlackBoxes

\vsize=7.5in

\def\pf{\hfill $\square$}
\def\c{\cite}

\def\fg{\frak{g}}

\def\end{\text{End}}

\def\ir{\int_{\Bbb R}}
\def\bk{\boldkey K}

\def\ba{\boldkey A}
\def\bb{\boldkey B}
\def\bc{\boldkey C}
\def\bi{\boldkey I}

\def\bs{\boldkey S}

\def\g{\gamma}
\def\lp{\ell_{+}}
\def\lm{\ell_{-}}

\topmatter
\title Factorization Problem on the Hilbert-Schmidt Group and
the Camassa-Holm Equation\endtitle
\leftheadtext{L.-C. Li}
\rightheadtext{Camassa-Holm Equation}

\author Luen-Chau Li\endauthor
\address{L.-C. Li, Department of Mathematics,Pennsylvania State University,
University Park, PA  16802, USA}\endaddress
\email luenli\@math.psu.edu\endemail

\abstract  In this paper, we solve the Camassa-Holm equation 
for a relatively large class of initial data by using a
factorization problem on the Hilbert-Schmidt group.
\endabstract
\endtopmatter

\document
\subhead
1. \ Introduction
\endsubhead

\baselineskip 15pt
\bigskip

The Camassa-Holm (CH) equation is a model of
long waves in shallow water.  Since its introduction in 
1993 by Camassa and Holm \c{CH}, the equation has received
considerable attention and its various aspects were
studied using a variety of methods. (See, for example, \c{BF},\c{BSS},
\c{Con}, \c{CM}, \c{CS}, \c{GH}, \c{M}, \c{XZ} and the references 
therein.)  

In contrast to the KdV equation, the CH equation admits breaking solutions.
However, a relatively large class of initial data which give rise to 
global solutions has also been identified in \c{C1} and independently in
\c{Con}.  In the  paper \c{C1} and its more elaborate version \c{C2}, 
the initial value problem for the CH equation was analyzed through
its characteristic formulation.  It is in this context that an isospectral 
problem in the 
form of an integro-differential equation was discovered for the
Lagrangian version of the integrable PDE. 
This remarkable fact has led, in particular, to a particle
method for numerically solving the CH equation.(See \c{CHL}
for subsequent development of this particle method.)

The main goal of this work is to show that the integro-differential
equation in the afore-mentioned papers of Camassa is in fact exactly solvable,
in the sense that a formula for its solution can be written down.
The upshot of this, as the reader will see, is that we can
integrate both  the Lagrangian version and the Eulerian version
of the CH equation explicitly.  The solution of the integro-
differential equation is based on a factorization problem
on the Hilbert-Schmidt group which we will introduce.
On the other hand, we will show that the factorization problem 
on the Hilbert-Schmidt group can be reduced to solving a
family of Fredholm integral equations, and this can be
achieved by using regularized Fredholm determinants
and  Fredholm first minors.  The reader will
see that the Lax equation which corresponds to
the integro-differential equation is in some
sense an infinite dimensional analog of the Toda
flow on $n\times n$ matrices (cf. \c{DLT}).  That such
an analog is connected with the Camassa-Holm equation
is rather remarkable and is responsible for
the elementary method of solving the equation here. 

The paper is organized as follows.  In Section 2, we begin by 
introducing a class of integrable isospectral deformations of 
Hilbert-Schmidt operators on $L^{2}(\Bbb R)$ using the r-matrix approach,
then we discuss the underlying Lie groups and coadjoint orbits.
Since the Hilbert-Schmidt operators on $L^{2}(\Bbb R)$ are given
by integral operators with kernels in $L^{2}(\Bbb R^{2})$, 
this leads naturally to a class of integrable integro-differential
equations.  In particular, for a special choice of Hamiltonian,
the corresponding Lax equation gives rise to the integro-differential
equation  obtained in \c{C1}, \c{C2}.  In Section 3, we discuss the solution 
of the factorization
problem on the Hilbert-Schmidt group and its application towards
the integration of the integro-differential equation.  Finally in 
Section 4, we show how to apply our 
results to obtain the explicit integration of the Lagrangian version and 
the Eulerian version of the Camassa-Holm equation. 

To close this introduction, we remark that the Lagrangian version
of the CH equation was also explicitly integrated in \c{M}. 
We note, however, that the assumption on the initial data
and the method used in \c{M} are quite different from the
one employed here.(See Remark 4.3 for the relationship
between the spectral problems.)

\medskip
\noindent{\bf Acknowledgments.} The main bulk of this work was completed 
at the Institut Mittag-Leffler when the author was participating
in the program on Wave Motion in the Fall of 2005. The author would 
like to thank the organizers for creating a wonderful atmosphere
during his visit and he acknowledges the financial support and
hospitality of the Institute.  He is especially indebted
to Adrian Constantin for the invitation and for bringing
the reference \c{Con} to his attention.  He also wishes to thank Boris
Kolev, Xavier Raynaud, and in particular Henry McKean for many
interesting discussions.  Finally, the author is most grateful to
the referee's helpful advice which has led to 
simplification of several arguments in Section 4.  He
also acknowledges the referee's helpful suggestions to
make the paper more readable and for the reference \c{C1}.

\bigskip
\bigskip

\subhead
2. \ Classical r-matrices and integrable integro-differential equations
\endsubhead

\bigskip

In this section, we introduce a class of integrable integro-differential
equations associated with isospectral deformations of Hilbert-Schmidt
operators.

Let ${\Cal H}$ be the Hilbert space $L^{2}({\Bbb R})$ consisting of
real-valued measurable functions on $\Bbb R$ that satisfy
$\ir |f(x)|^{2}\,dx < \infty$ with inner product
$$(f,g) = \ir f(x)g(x) dx \eqno(2.1)$$
and let  $\fg$ be the space of  Hilbert-Schmidt operators on ${\Cal H}$.
It is well-known that \c{RS} $\bk\in \fg$ if and only if there is
a function $K\in L^{2}(\Bbb R^2)$ uniquely determined by $\bk$ such that
$$(\bk \varphi)(x) = \ir K(x,y)\varphi(y)\,dy .\eqno(2.2)$$
In other words, the Hilbert-Schmidt operators on ${\Cal H}$ are
precisely the integral operators on ${\Cal H}$ with $L^{2}$-kernels.
Moreover, for $\bk\in \fg$, the Hilbert-Schmidt norm is given by
$$\|\bk \|^{2}_{2} = \int_{\Bbb R^{2}} |K(x,y)|^{2}\, dxdy.\eqno(2.3)$$

Let $B({\Cal H})$ be the space of bounded operators on ${\Cal H}$.
We recall that if $T\in B({\Cal H}),$ its adjoint $T^*$ is defined
by means of the relation
$$(T^{*}\varphi, \psi) = (\varphi, T\psi)\eqno(2.4)$$
for all $\varphi\in {\Cal H}$ and $\psi\in {\Cal H}.$
Specializing to the case where $\ba\in \fg \subset B({\Cal H})$ with
kernel $A(x,y)$, it is 
easy to show from (2.4) and
Fubini's theorem that its adjoint $\ba^{*}$ is the integral operator on 
$\Cal H$ with kernel $A^{*}(x,y)= A(y,x)$.  Thus $\ba^*$ is also in $\fg.$
In what follows, if $\ba \in \fg$, we shall use the notation
$$\ba \doteq A(x,y)\eqno(2.5)$$
to mean that the integral operator $\ba$ has kernel $A(x,y).$

\proclaim
{Proposition 2.1} $\fg$ is a Hilbert Lie algebra with Lie bracket
$[\cdot,\cdot]$ defined by
$$\eqalign{
[\ba, \bb] &= \ba\circ \bb - \bb\circ\ba \cr
           &\doteq \int_{\Bbb R} (A(x,z)B(z,y)-B(x,z)A(z,y))\,dz\cr}
\eqno(2.6)$$
and the inner product on $\fg$ is the usual Hilbert-Schmidt
inner product $(\cdot,\cdot)_{2}$, i.e.,
$$\eqalign{
(\ba, \bb)_{2}& = tr (\ba^{*}\circ\bb)\cr
             & =\int_{\Bbb R^{2}} A(x,y)B(x,y)\,dxdy.\cr}\eqno(2.7)$$
Moreover, $\fg$ is equipped with the non-degenerate
ad-invariant pairing
$$(\ba, \bb)= \int_{\Bbb R^{2}} A(x,y)B(y,x)\, dx dy.\eqno(2.8)$$
\endproclaim

\demo
{Proof} Since the Hilbert-Schmidt operators is closed
under the operations of addition, subtraction, and 
composition, it follows that
the bracket operation in (2.6) is well-defined and
it is clear that $[\cdot,\cdot]$ is a Lie bracket.
On the other hand, it is well-known that $\fg$ with
the inner product in (2.7) is a Hilbert space.
Hence $\fg$ is a Hilbert Lie algebra.  We shall
leave the rest of the assertion to the reader as
an exercise.
\pf
\enddemo  

In addition to what we have above, we remark that as a special
case of a  general theorem (see, for example, \c{RS}), $\fg$ 
is a 2-sided ideal in $B({\Cal H})$.
Now, let $\bi\in B({\Cal H})$ be the identity operator, and let $GL({\Cal H})$
denote the group of invertible operators in $B({\Cal H})$. 
We define
$$ G = GL({\Cal H})\bigcap\, (\bi + \fg).\eqno(2.9)$$
If  $\bi + \bk\in G$, it is well-known that $(\bi + \bk)^{-1}$
is also in $G$. (See, for example, \c{Sm}.) Hence $G$ is a group
under the composition of operators.   As a matter of fact, $G$
is a Hilbert Lie group which
integrates the Lie algebra $\fg,$ the Hilbert manifold
structure is being determined by the map
$G\longrightarrow \fg: \bi + \bk\mapsto \bk$
which is a bijection onto an open subset of $\fg$
consisting of Hilbert-Schmidt operators for which
$-1$ is not an eigenvalue.
We will call $G$
the {\it Hilbert-Schmidt group}.
In this case, the adjoint action of the group $G$  on $\fg$ is given
by the formula  
$$Ad_G(g)\bk =g\circ\bk\circ g^{-1}.\eqno(2.10)$$
Because the pairing on $\fg$ is ad-invariant, we also have
$$ (g\circ\ba \circ g^{-1}, \bb)=(\ba, g^{-1}\circ\bb\circ g).\eqno(2.11)$$
On the other hand, the exponential map $\exp:\fg\longrightarrow G$ is given 
by the expression
$$\exp(\bk) = \sum_{j=0}^{\infty} \frac{\bk^{j}}{j!}\eqno(2.12)$$
where the powers of $\bk$ are defined recursively by
$$\bk^{0} = \bi,\quad \bk^{j} = \bk\circ \bk^{j-1},\,\,j=1,2,\cdots .\eqno(2.13)$$

The Hilbert Lie algebra $\fg$ has two distinguished Lie subalgebras
$\frak l$ and $\frak k$, where $\frak l$ consists of Volterra integral
operators $\ba$ of the form
$$(\ba \varphi)(x) = \int_{-\infty}^{x} A(x,y)\varphi(y)\,dy\eqno(2.14)$$
and $\frak k$ consists of integral operators $\bb$ for which
$$\bb^{*} = - \bb.\eqno(2.15)$$
We will call $\frak l$ the {\it lower triangular subalgebra} of $\fg$
and $\frak k$ the {\it skew-symmetric subalgebra}.

\proclaim
{Proposition 2.2} We have
$$\fg = \frak l \oplus \frak k.$$
\endproclaim

\demo
{Proof} Given $\bk\in \fg$ with kernel $K(x,y)$, it is clear that
$$\bk = \Pi_{\frak k} \bk  + \Pi_{\frak l} \bk$$
where
$$\Pi_{\frak k} \bk \doteq K(x,y)\chi_{(x,\infty)}(y) - K(y,x)\chi_{(-\infty,x)}(y)$$
and 
$$\Pi_{\frak l} \bk \doteq K(x,y)\chi_{(-\infty,x)}(y) + K(y,x)\chi_{(-\infty,x)}(y).$$
In the above expressions, $\chi_{(x,\infty)}(y)$ and $\chi_{(-\infty, x)}(y)$ 
are respectively the characteristic
functions of $(x,\infty)$ and $(-\infty, x)$.
Now it is straightforward to check that
$\Pi_{\frak l} \bk \in \frak l$ and
$\Pi_{\frak k} \bk \in \frak k$.  Therefore, $\fg = \frak l + \frak k.$
To show that the sum is direct, suppose $K(\cdot,\cdot)$ is
the kernel of an operator which belongs to both $\frak l$ and
$\frak k$.  Then $K(x,y)=0$ away from the diagonal and hence
the corresponding operator is zero.  This completes the
proof.
\pf
\enddemo

Let $\Pi_{\frak l}$ and $\Pi_{\frak k}$ be the projection operators
onto $\frak l$ and $\frak k$ respectively associated with
the decomposition $\fg = \frak l\oplus \frak k.$  Then
it follows from \c{STS} that
$$ R = \Pi_{\frak l} - \Pi_{\frak k}\eqno(2.16)$$
is a classical r-matrix on $\fg$ satisfying the modified 
Yang-Baxter equation (mYBE) 
$$[R(\ba),R(\bb)] -R([R(\ba),\bb] + [\ba, R(\bb)]) = -[\ba,\bb]\eqno(2.17)$$
for all $\ba$, $\bb\in \fg.$  Consequently, the formula
$$[\ba, \bb]_{R} = {1\over 2}([R(\ba),\bb] + [\ba, R(\bb)]),\quad \ba,\bb\in
\fg \eqno(2.18)$$
defines a second Lie bracket on $\fg$  and we shall denote the associated
Lie algebra by  $\fg_{R}.$   In what follows, we shall compute
the dual maps of all linear operators on $\fg$ with respect
to the pairing $(\cdot,\cdot)$ in (2.8).

\proclaim
{Proposition 2.3} If $\ba \in \fg$ and $\ba \doteq A(x,y)$, then
$$\aligned
&\Pi^{*}_{\frak k} \ba \doteq (A(x,y)-A(y,x))\chi_{(-\infty,x)}(y)\\
&\Pi^{*}_{\frak l} \ba \doteq A(x,y)\chi_{(x,\infty)}(y) +
 A(y,x)\chi_{(-\infty,x)}(y).\\
\endaligned
$$
\endproclaim

The proof is a straightforward calculation and so we skip the
details.  We shall equip $\fg^{*}_{R}\simeq \fg$ with the Lie-Poisson
structure
$$\{\,F_1,F_2\,\}_{R}(\bk) =(\bk, [dF_1(\bk),dF_2(\bk)]_{R})\eqno(2.19)$$
where $F_1$, $F_2 \in C^{\infty}(\fg^{*}_{R}),$ and
$dF_{i}(\bk)\in \fg$ is defined by the formula
${d\over dt}\big|_{t=0} F_{i}(\bk + t\bk^{\prime}) = (dF_{i}(\bk), \bk^{\prime}),\,
i=1,2.$

The following result is a consequence of standard classical r-matrix theory. 
(See  \c{STS} and \c{RSTS} for the general theory.)

\proclaim
{Proposition 2.4} (a) The Hamiltonian equations of motion generated
by $F \in C^{\infty}(\fg^{*}_{R})$ is given by
$$\dot \bk = {1\over 2}[R(dF(\bk)),\bk] -{1\over 2}R^{*}[\bk, dF(\bk)].
\eqno(2.20)$$
In particular, for the Hamiltonian $H_{j}(\bk) = {1\over 2(j+1) }tr(\bk^{j+1})$,
$j= 1,2,\ldots$, the corresponding equation is the Lax equation
$$\dot \bk = {1\over 2} [\,\Pi_{\frak l} \,\bk^{j},\bk\,].\eqno(2.21)$$
\newline
(b) The family of functions $H_{j}(\bk)$, $j=1,2,\ldots$ Poisson
commute with respect to $\{\cdot,\cdot\}_{R}$.
\endproclaim

\demo
{Proof} The Hamiltonian equation of motion (2.20) is obtained
from (2.19) by a direct calculation.  On the other hand, by 
using (2.11), we find that
$Ad^{*}_{G} (g^{-1})\ba = g\circ\ba\circ g^{-1}$, $g\in G.$   Since
$tr (\ba\circ\bb) = tr (\bb\circ\ba)$ for any
$\bb\in B({\Cal H})$ and any trace class operator $\ba$,
it follows that $H_{j}(Ad^{*}_{G}(g^{-1})\bk) = H_{j}(\bk),$ $g\in G.$
By classical r-matrix theory, we then conclude that
the family of functions  $H_{j}(\bk)$, $j=1,2,\ldots$
Poisson commute with respect to $\{\cdot,\cdot\}_{R}$.
The equation of motion for $H_j$ now follows from (2.20)
as the invariance property  $H_{j}(Ad^{*}_{G}(g^{-1})\bk) = H_{j}(\bk)$
implies  $[\bk, dH_{j}(\bk)]=0.$  This completes the proof.
\pf
\enddemo
\medskip
Let
$$\frak p = \lbrace \bk\in \fg\mid \bk = \bk^{*} \rbrace.\eqno(2.22)$$

\proclaim
{Corollary 2.5} (a) $\frak p$ is a Poisson submanifold of
$(\fg^{*}_{R}, \{\cdot,\cdot\}_{R})$. Hence eqn. (2.21) with
$\bk\in \frak p$ is Hamiltonian with respect to the
induced Poisson structure on $\frak p.$
\newline
(b) For the Hamiltonian $H_{1}(\bk) = {1\over 4}(\bk,\bk)$,
the evolution of the kernel $K(x,y;t)$ corresponding to $\bk$ is
given by the integro-differential equation
$$\eqalign{
{\dot K}(x,y\,;t)  = & {1\over 2}\int_{-\infty}^{x} K(x,z\,;t)K(z,y\,;t)\,dz -
 {1\over 2} \int_{y}^{\infty} K(x,z\,;t)K(z,y\,;t)\,dz\cr
 & + {1\over 2}\int_{-\infty}^{x} K(z,x\,;t)K(z,y\,;t)\,dz -
   {1\over 2}\int_{y}^{\infty} K(x,z\,;t)K(y,z\,;t) \,dz.\cr}\eqno(2.23)$$
In the special case when $\bk$ belongs to the Poisson submanifold
$\frak p$, the corresponding kernel $K(x,y;t)$ is symmetric. In
this case, the above equation reduces to
$$\eqalign{
{\dot K}(x,y\,;t) & =\int_{-\infty}^{x} K(x,z\,;t)K(z,y\,;t)\,dz -
\int_{y}^{\infty} K(x,z\,;t)K(z,y\,;t)\, dz\cr
             & =\int_{-\infty}^{y} K(x,z\,;t)K(z,y\,;t)\,dz -\int_{x}^{\infty}
                K(x,z\,;t)K(z,y\,;t)\, dz.\cr}\eqno(2.24)$$
\endproclaim

\demo
{Proof} (a) From (2.20), the Hamiltonian vector field generated
by $F$ can be rewritten in the form
$$X_{F}(\bk) = [\bk, \Pi_{\frak k}(dF(\bk))] -\Pi^{*}_{\frak l}[\bk, dF(\bk)].$$
From the expression for $\Pi^{*}_{\frak l}$ in Proposition 2.3, it
is clear that the second term in the above expression is always
in $\frak p$.  On the other hand, if $\bk \in \frak p$, then
it is easy to check that $[\bk, \Pi_{\frak k}(dF(\bk))]$ is also
in $\frak p.$  Consequently, we have $X_{F}(\bk)\in \frak p$
for $\bk\in \frak p$ and this shows that $\frak p$ is a
Poisson submanifold of $(\fg^{*}_{R}, \{\cdot,\cdot\}_{R})$.
\newline
(b) We have
$$(\Pi_{\frak l}\bk\circ \bk)\varphi(x) =\int_{\Bbb R}\left(\int_{-\infty}^{x}
   (K(x,z\,;t)+ K(z,x\,;t))K(z,y\,;t)\,dz\right) \varphi (y)\, dy$$
while
$$(\bk\circ \Pi_{\frak l}\bk)\varphi(x) =\int_{\Bbb R}\left(\int_{\Bbb R}
   K(x,z\,;t)(K(z,y\,;t)+K(y,z\,;t)) \chi_{(-\infty, z)}(y)\,dz\right) 
  \varphi (y)\,dy.$$
Therefore, the evolution of the kernel $K(x,y\,;t)$ is given by
$$\aligned
2{\dot K}(x,y\,;t)  = & \int_{-\infty}^{x} K(x,z\,;t)K(z,y\,;t)\,dz -
  \int_{y}^{\infty} K(x,z\,;t)K(z,y\,;t)\,dz\\
 & + \int_{-\infty}^{x} K(z,x\,;t)K(z,y\,;t)\,dz -
   \int_{y}^{\infty} K(x,z\,;t)K(y,z\,;t)\,dz.\\
\endaligned
$$
\pf
\enddemo

For our next remark and the discussion in Section 3, we  
introduce the Lie subalgebra $\frak u$ of $\fg$ which consists of
Volterra integral operators $\bb$ of the form
$$(\bb \varphi)(x) = \int_{x}^{\infty} B(x,y)\varphi(y)\,dy. \eqno(2.25)$$
\medskip
\noindent{\bf Remark 2.6.} (a) From the definition of $\frak k$ and
$\frak l$, and from equation (2.21), it is clear that what 
we are dealing with here is in some sense an infinite dimensional
analog of the Toda flows on $n\times n$ matrices (cf.~ \c{DLT}).
\smallskip
\noindent (b) A different decomposition of the Hilbert Lie
algebra $\fg$ is given by 
$$\fg = \frak l \oplus \frak u\eqno(2.26)$$
with associated projection maps $\Pi_{-}:\fg\longrightarrow \frak l$
and $\Pi_{+}:\fg\longrightarrow \frak u.$ (We will use these
projection maps in Section 3.)
Indeed, if we consider the r-matrix $R= \Pi_{-} -\Pi_{+}$
associated with this splitting and equip the corresponding
$\fg^{*}_{R}$ with the Lie-Poisson structure, then the
evolution of the {\it general kernel} $K(x,y;t)$ under the Hamiltonian
flow generated by ${1\over 2}(\bk, \bk)$ is given by
the equation
$$\aligned
{\dot K}(x,y\,;t) & =\int_{-\infty}^{x} K(x,z\,;t)K(z,y\,;t)\,dz -
\int_{y}^{\infty} K(x,z\,;t)K(z,y\,;t)\, dz\\
             & =\int_{-\infty}^{y} K(x,z\,;t)K(z,y\,;t)\,dz -\int_{x}^{\infty}
                K(x,z\,;t)K(z,y\,;t)\, dz.
\endaligned
$$
Clearly, eqn.(2.24) is a special case of this.  Note, 
however, that $\frak p$ is not a Poisson submanifold any more
with this choice of r-matrix and the corresponding Lie-Poisson
structure. Indeed, the Hamiltonian vector field generated
by a {\it general function} $F$ is now of the form 
$[\bk, \Pi_{+}(dF(\bk))] -\Pi^{*}_{-}[\bk, dF(\bk)].$
Clearly, this is not necessarily in $\frak p$ for
$\bk \in \frak p.$
Thus from the Hamiltonian point of view, the r-matrix in (2.16)
is the correct choice. For the relation
between (2.24) and the Camassa-Holm equation, and the
coadjoint orbit picture, we refer the reader to the
discussion in Section 4 preceding Remark 4.1, Remark 4.2 and 
Proposition 2.8 below for details.
\medskip

In the rest of the section, we shall describe the symplectic
leaves of the Lie-Poisson structure $\{\cdot,\cdot \}_{R}$ which
are given by the coadjoint orbits of the infinite dimensional
Lie group $G_{R}$  which integrates $\fg_R.$   In particular,
we shall consider the coadjoint action of $G_R$ on the class $\frak p_{*}$ of
Hilbert-Schmidt operators $\bk\in \fg$ with so-called
single-pair kernels \c{GK}.  By definition, 
a Hilbert-Schmidt operator $\bk\in \frak p_{*}$ if and
only if its kernel is of the form
$$K(x,y) =\cases a(x)b(y), & x\leq y \\
                 a(y)b(x), & x>y, \endcases\eqno(2.27)$$
where $a$ and $b$ are functions on $\Bbb R.$ (Note that
$a$ and $b$ are not necessarily in $L^{2}(\Bbb R)$.)  In order to 
describe $G_R$, we begin by introducing the Lie subgroups (see
Remark 2.7)
$$\eqalign{
&{\Cal L} = \bi + \frak l, \cr
&{\Cal K}= \{\, k\in G\mid k\circ k^{*} = k^{*}\circ k = \bi\, \}\cr}
\eqno(2.28)$$
of $G$ which corresponds to the Lie algebras $\frak l$ and
$\frak k.$ 
\medskip
\noindent{\bf Remark 2.7.}  It is clear that the group operation
of $G$ is closed on ${\Cal L}$.  On the other hand, if $\ba\in \frak l,$
we can show that the Neumann series $\sum_{j=0}^{\infty} (-1)^{j}\ba^{j}$
converges by using the estimate 
$\|\ba^{j+1}\|\leq {\|\ba\|^{j+1}_{2}}/{\sqrt{(j-1)!}},$ $j=1,2,\ldots,$
which we can derive from eqn. (8), Section 2.7 of \c{Sm}
provided we interpret the inequality there as being valid almost
everywhere under our weaker assumption.   Thus $\bi + \ba$ is
invertible and 
$(\bi + \ba)^{-1} = \sum_{j=0}^{\infty} (-1)^{j}\ba^{j}\in {\Cal L}.$   This 
shows that
${\Cal L}$ is a subgroup of $G.$  That ${\Cal L}$ is a
Lie subgroup of $G$ now follows since the former is clearly
a submanifold of the latter.
\medskip
Let
$$G_{R} = \lbrace g\in G\mid g = g_{-}\circ g^{-1}_{+}, \,\,\hbox{where}\,\,
  g_{-}\in {\Cal L}, \, g_{+}\in {\Cal K}\, \rbrace.\eqno(2.29)$$
Then following the procedure in \c{DLT}, we can endow $G_{R}$ with a Lie 
group structure
by defining the multiplication
$$g\ast h \equiv g_{-}\circ h\circ g_{+}^{-1}\eqno(2.30)$$
and we can show that $(G_{R}, \ast)$ is a Lie group which
corresponds to the Lie algebra $\fg_{R}.$  Moreover,
the adjoint action of $G_{R}$ on $\fg_{R}$ is given by
$$Ad_{G_{R}}(g)\bk = g_{-}\circ\Pi_{\frak l}\bk\circ g^{-1}_{-}
  + g_{+}\circ\Pi_{\frak k}\bk\circ g^{-1}_{+}.\eqno(2.31)$$
Hence an easy computation using (2.11) shows that
$$Ad^{*}_{G_{R}}(g^{-1})\bk = \Pi^{*}_{\frak l}(g_{-}\circ \bk\circ g^{-1}_{-})
  + \Pi^{*}_{\frak k}(g_{+}\circ \bk\circ g^{-1}_{+})\eqno(2.32)$$
and the symplectic leaves of $\{\cdot,\cdot \}_{R}$ are
the orbits of this coadjoint action.

\proclaim
{Proposition 2.8} The class $\frak p_{*}\subset \frak p$ of
Hilbert-Schmidt operators with single-pair kernels
is invariant under $Ad^{*}_{G_R}.$
\endproclaim
 
\demo
{Proof} Take $\bk\in \frak p_{*}$,
$$\bk\doteq K(x,y)=\cases a(x)b(y), & x\leq y \\
                    a(y)b(x), & x>y. \endcases$$
Then for $g = g_{-}\circ g^{-1}_{+}\in G_{R}$, it is clear that
$g^{-1}_{+}\circ \bk\circ g_{+}\in \frak p.$  Therefore,
$\Pi^{*}_{\frak k} (g^{-1}_{+}\circ \bk\circ g_{+})=0$ so
that
$$Ad^{*}_{G_{R}}(g)\bk = \Pi^{*}_{\frak l}(g^{-1}_{-}\circ \bk\circ g_{-}).$$
Now, by a straightforward computation using the form
of $K(x,y)$ above and the fact that $g_{-}\in {\Cal L}$, we can show 
$$\Pi^{*}_{\frak l}(g^{-1}_{-}\circ \bk\circ g_{-})\doteq
\cases (g^{-1}_{-}a)(x)(g^{*}_{-}b)(y), & x\leq y \\
       (g^{-1}_{-}a)(y)(g^{*}_{-}b)(x), & x>y \endcases$$
from which we conclude that 
$Ad^{*}_{G_{R}}(g)\bk\in \frak p_{*},$ as asserted.
\pf
\enddemo

From this result, it follows that the coadjoint orbit of $G_{R}$ through an
element $\bk\in \frak p_{*}$ will consist entirely of elements from 
$\frak p_{*}.$
In particular, this means that if the initial data of (2.24)
is a single-pair kernel, then $K(x,y;t)$ is also
a single-pair kernel for all $t.$   This is the fact 
which underlies
the geometry behind our application in Section 4 below.
\smallskip
\noindent{\bf Remark 2.9.} (a) As the reader will see
in Theorem 3.1, every $g\in G$ admits a unique factorization
$g= g_{-}\circ g^{-1}_{+}$ with $g_{-}\in {\Cal L}$ and 
$g_{+}\in {\Cal K}.$  Thus the underlying manifold
of $(G_{R}, \ast)$ is just $G$ in this case.
Note that this factorization result in Theorem 3.1 is
nothing but the global version of the decomposition
in Proposition 2.2.
\smallskip
\noindent (b) As a final remark of this section,
we would like to point out that everything we have done
in this section can be pushed through in the more general setting
when the Hilbert space is taken to be $L^{2}(\Bbb R, d\mu)$
for an arbitrary Borel measure $\mu$ on $\Bbb R.$

\bigskip
\bigskip

\subhead
3. Solution by factorization
\endsubhead
\bigskip

We recall the Lie subgroups

$$\eqalign{
&{\Cal L} = \bi + \frak l, \cr
&{\Cal K}= \{\, k\in G\mid k\circ k^{*} = k^{*}\circ k = \bi\, \}\cr}
\eqno(3.1)$$
of $G$ introduced in Section 2.  As the reader will see, they
play an important role in the solution of (2.24).  In order to
discuss the factorization problem, let us also recall several formulas from 
the theory of
regularized determinants which we are going to need in
our context.  The reader is referred to \c{S}, \c{Sm} for
more details.

Let $\ba$ be a Hilbert-Schmidt operator on ${\Cal H} = L^{2}(\Bbb R)$,
then 
$${\Cal R}_{2}(\ba) := (\bi + \ba) e^{-\ba} - \bi\eqno(3.2)$$ 
is of trace class.  Following \c{S}, we can define the regularized determinant
$$\eqalign{
{\det}_{2}(\bi + \ba)& := \det (\bi + {\Cal R}_{2}(\ba))\cr
                 & := \sum_{k=0}^{\infty} tr \wedge^{k}({\Cal R}_{2}(\ba))\cr}
\eqno(3.3)$$
which obeys the estimate
$$|{\det}_{2}(\bi + \ba)| \leq \exp(\|{\Cal R}_{2}(\ba)\|_{1})  \eqno(3.4)$$
where $\|{\Cal R}_{2}(\ba)\|_{1} = 
tr (\sqrt {{\Cal R}_{2}(\ba)^*{\Cal R}_{2}(\ba)}\,).$
Suppose $\bi + \bk \in G$, then the analog of the first Fredholm minor
is given by
$$D_{2}(\bk) = -\bk(\bi + \bk)^{-1} {\det}_{2}(\bi + \bk)\eqno(3.5)$$
and so 
we have the formula
$$(\bi + \bk)^{-1} = \bi + \frac{D_{2}(\bk)}{{\det}_{2}(\bi + \bk)}.\eqno(3.6)$$
In connection with (3.6) above, it is important to note the Plemelj-Smithies
formulas
$${\det}_{2}(\bi + \bk) = 1 + \sum_{m=1}^{\infty} \alpha^{(2)}_{m} (\bk)\eqno(3.7)$$
and
$$D_{2}(\bk) = \bk +  \sum_{m=1}^{\infty} \beta^{(2)}_{m} (\bk).\eqno(3.8)$$
Here, 
$$ \alpha^{(2)}_{m} (\bk) = {1\over m!} \,\det \pmatrix
   0 & m-1 & 0 &\cdots & 0\cr
   \sigma_{2}(\bk) & 0 & m-2 &\cdots & 0\cr
   \sigma_{3}(\bk) &\sigma_{2}(\bk)&0 &\cdots & 0\cr
   \vdots & \vdots & \vdots & \cdots & \vdots \cr
   \sigma_{m-1}(\bk) & \sigma_{m-2}(\bk) & \sigma_{m-3}(\bk)& \cdots& 1\cr
   \sigma_{m}(\bk) & \sigma_{m-1}(\bk) & \sigma_{m-2}(\bk)&\cdots& 0 \cr
  \endpmatrix  \eqno(3.9)$$
and
$$\beta^{(2)}_{m} (\bk) = {1\over m!} \,\det \pmatrix
   \bk & m  & 0 & \cdots & 0\cr
   \bk^{2} & 0 & m-1 & \cdots & 0\cr
   \bk^{3} & \sigma_{2}(\bk) & 0 & \cdots & 0\cr
   \vdots & \vdots & \vdots &\cdots & \vdots \cr
   \bk^{m} & \sigma_{m-1}(\bk) & \sigma_{m-2}(\bk) & \cdots & 1\cr
   \bk^{m+1} & \sigma_{m}(\bk) & \sigma_{m-1}(\bk) & \cdots & 0\cr
   \endpmatrix \eqno(3.10)$$
where
$$\sigma_{j}(\bk) = tr (\bk^{j}), \quad j \geq 2.\eqno(3.11)$$

We next introduce a piece of notation.  If $\bk \in \fg$ and
$y\in \Bbb R$, we shall denote by $\bk \mid_{ (-\infty,y)}$
the operator $L^{2}(-\infty,y)\longrightarrow L^{2}(-\infty,y)$
defined by
$$(\bk\mid_{(-\infty,y)} \varphi)(x) = \int_{-\infty}^{y} K(x,z)\varphi(z)\, dz 
\eqno(3.12)$$
for $\varphi\in L^{2}(-\infty,y).$  Similarly, $(\bi + \bk)\mid_{(-\infty,y)}$
has an analogous meaning.

In order to solve the integro-differential equation (2.24), the
following result is basic.  

\proclaim
{Theorem 3.1} Suppose $\bi + \bk \in G,$
then $\bi + \bk$ has a unique factorization
$$\bi + \bk = b_{-}\circ b_{+}^{-1} \eqno(3.13)$$
where $b_{-}\in {\Cal L}$ and
$b_{+}  \in {\Cal K}.$  If $(b^{-1}_{-})^{*}-\bi \doteq C_{+}(x,y)$, 
then explicitly, $C_{+}(x,y)$ is given by
$$\eqalign{
C_{+}(x,y) & = -\left(((\bi + \bs)|_{(-\infty, y)})^{-1} S(\cdot,y)\right)(x),\quad
x<y \cr
& = -S(x,y) -\frac {(D_{2}(\bs|_{(-\infty, y)})S(\cdot,y))(x)}
    {{\det}_{2}((\bi + \bs)|_{(-\infty, y)})}\cr}\eqno(3.14)$$
and $C_{+}(x,y) = 0$ for $y<x$ where $S(x,y)$ is the kernel of
$$\bs = \bk + \bk^{*} + \bk\circ\bk^{*}.\eqno(3.15)$$
\endproclaim

\demo
{Proof} By the analog of Remark 2.7,
$${\Cal U} = \bi + \frak u$$
is the Lie group corresponding to the Lie algebra $\frak u$
introduced at the end of Section 2.
For each $y\in \Bbb R$, consider the equation
$$C(x,y) + \int_{-\infty}^{y} S(x,z)C(z,y)\,dz = -S(x,y),\quad  x<y.
\eqno(3.16)$$
Since $\bi + \bs = (\bi + \bk)\circ (\bi + \bk)^{*}$ is positive
definite, it follows that $(\bi + \bs)|_{(-\infty, y)}$ is invertible.
Hence (3.16) has a unique solution given by
$$\aligned
C_{+}(x,y) & = -\left(((\bi + \bs)|_{(-\infty, y)})^{-1} S(\cdot,y)\right)(x),\quad
x<y \\
& = -S(x,y) -\frac {(D_{2}(\bs|_{(-\infty, y)})S(\cdot,y))(x)}
    {{\det}_{2}((\bi + \bs)|_{(-\infty, y)})}\\
\endaligned$$
where we have used the formula in (3.6). Set $C_{+}(x,y) = 0$ for
$y<x$ and  let $\bc_{+}$ denote
the corresponding operator in $\frak u$.  Then from (3.16), we have
$$\bc_{+} +\,\Pi_{+}(\bs\circ \bc_{+}) = -\bs_{+} \eqno(3.17)$$
where $\Pi_{+}:\fg\longrightarrow \frak u$ is the projection
operator to $\frak u$ relative to the splitting in (2.26)
and $\bs_{+} = \Pi_{+} \bs.$
But on the other hand, we find
$$\eqalign{
& (1-\Pi_{+})(\bs + \bs\circ \bc_{+}) - \bs\cr
= & \,(\bi + \bs)\circ \bc_{+}.\cr}\eqno(3.18)$$
Hence it follows that
$$(\bi + \bs)\circ (\bi + \bc_{+}) =\, \bi + \bb_{-}$$
where 
$$\bb_{-}:=\Pi_{-}(\bs +\bs\circ\bc_{+})\in \frak l.\eqno(3.19)$$
Set $b_{-} = \bi + \bb_{-}$ and $c_{+} = \bi + \bc_{+}$.
Then $\bi + \bs = b_{-}\circ c^{-1}_{+}$.  But from 
$(\bi + \bs)^{*}= \bi + \bs$, we also have
$\bi + \bs = (c^{-1}_{+})^{*}\circ b^{*}_{-}.$  Equating
the two expression for $\bi + \bs$, we find
$$c^{*}_{+}\circ b_{-} = b^{*}_{-}\circ c_{+}\in {\Cal L}\cap {\Cal U}.$$
As  ${\Cal L}\cap {\Cal U} = {\bi}$,  we conclude that
$c_{+} = (b^{*}_{-})^{-1}$ and so we have established the factorization
$$\bi + \bs = b_{-}\circ b^{*}_{-}.$$  
Now we define
$$b_{+} = (\bi + \bk)^{-1}\circ b_{-}.$$
Then a straightforward verification shows that
$b_{+}\in {\Cal K}.$  Finally, the uniqueness of the factors
$b_{\pm}$ is obvious.
\pf
\enddemo

\proclaim
{Theorem 3.2} Let $\bk_0\in \fg$ and let 
$b_{-}(t)\in {\Cal L}$, $b_{+}(t)\in {\Cal K}$ be the unique solution of the
factorization problem
$$\exp\left(-{1\over 2}t\bk_0\right) = b_{-}(t)\circ b_{+}(t)^{-1}.\eqno(3.20)$$
Then for all $t$,
$$\bk(t) = b_{\pm}(t)^{-1}\circ \bk_{0}\circ b_{\pm}(t) \eqno(3.21)$$
solves the initial value problem
$$\dot \bk = {1\over 2}[\,\Pi_{\frak l}\,\bk,\bk\,]= 
{1\over 2}[\,\bk, \Pi_{\frak k}\,\bk], \quad 
\bk(0) = \bk_0.
\eqno(3.22)$$
\endproclaim

\demo
{Proof}  We shall present a direct proof of the theorem.  First of all,
the factorization problem in (3.20) has unique solutions $b_{-}(t)\in {\Cal L}$
and $b_{+}(t)\in {\Cal K}$ by Theorem 3.1. Take
$$\bk(t) = b_{+}(t)^{-1}\circ \bk_{0}\circ b_{+}(t).$$
By differentiating the expression, we have
$$\dot \bk(t) = [\, \bk(t),\, b_{+}(t)^{-1}\circ \dot b_{+}(t)\,].\qquad \qquad
\qquad \qquad (*)$$
On the other hand, by differentiating (3.20), we find
$$-{1\over 2}\bk(t) = b_{-}(t)^{-1}\circ \dot b_{-}(t) 
- b_{+}(t)^{-1}\circ \dot b_{+}(t).$$
Hence by applying $\Pi_{\frak k}$ to both sides of the above expression,
the result is
$$b_{+}(t)^{-1}\circ \dot b_{+}(t) = {1\over 2}\Pi_{\frak k}\,\bk(t) .$$
Therefore, on substituting into (*), we conclude that
$$\dot \bk(t) = {1\over 2} [\, \bk(t), \Pi_{\frak k}\,\bk(t) \,].$$
This shows that $\bk(t) = b_{+}(t)^{-1}\circ \bk_{0}\circ b_{+}(t)$ solves
the initial value problem.\pf
\enddemo

\medskip
\noindent{\bf Remark 3.3.}  (a)  The factorization method is the
most important feature of classical r-matrix theory.  It should
be emphasized that there is no universal method to solve the
factorization problems.  Rather, the method varies with the
Lie groups involved. For examples involving finite dimensional
matrix groups, the reader is referred to \c{DLT} for further
information.  On the other hand, factorization problems
associated with loop groups are related to Riemann-Hilbert
problems. See, for example, \c{RSTS} and \c{DL} in this
connection.
\newline
(b) We can also give a geometric proof of Theorem 3.2 by using
Poisson reduction.  Indeed, from this point of view, we can
understand the Hamiltonian flow in (3.21) as the projection of
some simple Hamiltonian flow on the cotangent bundle $T^{*}G$.
We will give a sketch of the argument here, following essentially
the outline on p. 180 of \c{STS}.  Let $(G_{R}, \ast)$ be the
Lie group introduced in Section 2 with Lie algebra  $\fg_{R}.$   
Consider the action of $G_R$ on $G$,
$$G_{R}\times G\longrightarrow G: \,\,g\cdot h = g_{-}\circ h\circ g^{-1}_{+},
  \,\,g\in G_{R}, h\in G.\eqno(3.23)$$
We can lift this up to obtain a canonical action on $T^{*}G$ \c{AM}.
Indeed, this lifted action is given by
$$G_{R}\times T^{*}G\longrightarrow T^{*}G: g\cdot (h,\bk)=
  (g_{-}\circ h\circ g^{-1}_{+}, Ad^{*}_{G}(g^{-1}_{+})\bk)\eqno(3.24)$$ 
where we have made the identification 
$T^{*}G\simeq G\times \fg^{*} \simeq G\times \fg$
by using left translation and the pairing on $\fg.$  Consequently,
by Poisson reduction, the orbit space $T^{*}G/G_{R}\simeq \fg^*\simeq \fg$
has a unique Poisson structure such that the canonical projection
$$\pi:T^{*}G\longrightarrow T^{*}G/G_{R}\simeq \fg, (g,\bk)\mapsto
  Ad^{*}_{G}(g_{+})\bk = g^{-1}_{+}\circ \bk\circ g_{+}\eqno(3.25)$$
is a Poisson map.  By a direct computation, we can show that
the Poisson structure on $\fg\simeq T^{*}G/G_{R}$ is nothing 
but $\{\cdot,\cdot\}_{R}$.  Now we consider the bi-invariant Hamiltonian
on $T^{*}G\simeq G\times \fg,$ given by 
$$\widehat{H}_{1}(g,\bk) = H_{1}(\bk).\eqno(3.26)$$
Then its Hamiltonian equations of motion are
$$\eqalign{ \dot g &= -{1\over 2} g\circ \bk,\cr
             \dot \bk &= 0.\cr}\eqno(3.27)$$ 
Therefore, if we denote the corresponding flow by $F_{t}$, we have
in particular that
$$F_{t}(\bi,\bk_{0}) = \left(\exp\left(-{1\over 2}t\bk_{0}\right), \bk_{0}\right).
\eqno(3.28)$$
Consequently, the Hamiltonian flow generated by the reduced
Hamiltonian $H_{red.} = H_{1}$ on the orbit space $T^{*}G/G_{R}\simeq \fg$
is given by
$$\eqalign{\bar{F}_{t}(\bk_{0})&= \pi\circ F_{t}(\bi,\bk_{0})\cr
           &= Ad^{*}_{G}(b_{+}(t))\bk_{0}\cr
           &= b_{+}(t)^{-1}\circ \bk_{0}\circ b_{+}(t),\cr}\eqno(3.29)$$
as required.
\medskip

We are now going to write down the solution of the integro-differential
equation (2.24) by combining the above theorems.  For this purpose,
let $\bk_0\doteq K(x,y)$, $\bk(t) \doteq K(x,y;t)$ and write
$$\exp(-t\bk_0) = \bi + \bs(t), \quad \bs(t) \doteq S(x,y;t).\eqno(3.30)$$
By Theorem 3.1, the solution  $b_{-}(t)$ of the factorization
problem (3.20) is such that
$$(b_{-}(t)^{-1})^{*} - \bi = \bc^{t}_{+} \doteq C_{+}(x,y;t)\eqno(3.31)$$
where
$$\eqalign{
C_{+}(x,y;t) & = -\left((e^{-t\bk_{0}}|_{(-\infty, y)})^{-1} S(\cdot,y;t)\right)(x)
\,\chi_{(x,\infty)}(y)\quad\cr
& = -\left(S(x,y;t) +\frac {(D_{2}(\bs(t)|_{(-\infty, y)})S(\cdot,y;t))(x)}
    {{\det}_{2}\left(e^{-t\bk_{0}}|_{(-\infty, y)}\right)}\right)\chi_{(x,\infty)}(y).\cr}
\eqno(3.32)$$
On the other hand, it follows from the proof of Theorem 3.1 (see
(3.19) and the definition of $b_{-}$) that
$$b_{-}(t) = \bi + \Pi_{-}(\bs(t) +\bs(t)\circ\bc^{t}_{+})\eqno(3.33)$$
and hence
$$b_{-}(t) - \bi \doteq \left(S(x,y;t) +\int_{-\infty}^{\infty} 
S(x,z;t)C_{+}(z,y;t) dz\right)\, \chi_{(-\infty,x)}(y).\eqno(3.34)$$
Consequently, from
$$\eqalign{
\bk(t) & =  b_{-}(t)^{-1}\circ \bk_{0}\circ b_{-}(t)\cr
       & = \bk_0 + (b_{-}(t)^{-1} -\bi)\circ\bk_0 +
         \left(\bk_0 + (b_{-}(t)^{-1} -\bi)\circ\bk_0\right)\circ 
          (b_{-}(t) -\bi)\cr}\eqno(3.35)$$
and (3.31), (3.34), we find
$$\eqalign{
K(\xi,\eta;t) & = K(\xi,\eta) + \int_{-\infty}^{\infty} C_{+}(\zeta, \xi;t)
K(\zeta,\eta)\, d\zeta \cr
& \quad +\int_{\eta}^{\infty} \left(K(\xi,\zeta_{2}) +\int_{-\infty}^{\infty} 
C_{+}(\zeta_{1}, \xi;t) K(\zeta_1,\zeta_2)\,d\zeta_1\right)\cr
&\quad \quad \quad \quad\cdot\left(S(\zeta_2,\eta;t) +\int_{-\infty}^{\infty} 
S(\zeta_2,\zeta_3;t)
C_{+}(\zeta_3,\eta;t)\,d\zeta_3\right)\,d\zeta_2 .\cr}\eqno(3.36)
$$
\smallskip
\noindent{\bf Remark 3.4.}  In a similar fashion, we can write
down the solution of the integro-differential equation (2.23).
Indeed, all we have to do is to replace $e^{-t \bk_0}$ in
(3.30) and (3.32) by 
$\exp\left(-{1\over 2}t\bk_0\right)\circ\exp\left(-{1\over 2}t(\bk_0)^*\right).$
The operator $\bs(t)$ and its kernel $S(x,y;t)$ are of course
much more complicated in this case.

\bigskip
\bigskip

\subhead
4. Solution of the Camassa-Holm equation
\endsubhead

\bigskip

In this section, we shall consider the Camassa-Holm (CH) equation in
the non-dispersive case:
$$u_t -u_{xxt} + 3uu_x = 2u_{x}u_{xx} + uu_{xxx}.\eqno(4.1)$$

We will begin with a sketch of the connection between (4.1) and 
the integro-differential equation (2.24), as first discovered by
Camassa \c{C1}.  To do so, we introduce the auxiliary variable
$$m(x,t) = (1 - \partial^{2}_{x}) u(x,t) \eqno(4.2)$$
and rewrite the CH equation in the form
$$m_t + u m_x = -2 m u_x.\eqno(4.3)$$
We shall make the following assumption on the initial data:

\noindent (i) $u_{0} = u(\cdot, 0)$ is in the Schwarz class
${\Cal S}(\Bbb R)$,

\noindent (ii)  $m(x,0) = u_{0}(x) - u''_{0}(x) >0$ for all
$x\in \Bbb R.$
\smallskip
In the Lagrangian point of view of the CH equation, we consider
the trajectory $q(\xi,t)$ of a fluid particle which at $t=0$ is
located at $\xi \in \Bbb R$. 
Thus we have
$$\dot q(\xi, t) = u(q(\xi,t),t), \quad q(\xi,0) =\xi\eqno(4.4)$$
and a straightforward calculation (see \c{C1}, \c{Con})  
using (4.3) and (4.4) shows 
that 
$$m(q(\xi,t),t) (q_{\xi}(\xi, t))^{2} = m(\xi, 0).\eqno(4.5)$$
Now it follows from \c{C1},\c{C2} and \c{Con} that for the class of 
initial data introduced above, we have
$$0 < q_{\xi}(\xi,t) < \infty \eqno(4.6)$$
for all $t >0.$   
In particular, this means that the trajectories of the
fluid particles never cross.  

Therefore, if we define 
$$w(x,t) = \int_{-\infty}^{x} \sqrt{m(y,t)}\,dy , \eqno(4.7)$$
then 
$$w(q(\xi,t),t) = w(\xi,0) =w_{0}(\xi) \eqno(4.8)$$
so that we can rewrite (4.5) as
$$m(q(\xi,t),t) = \left(\frac{w^{\prime}_{0}(\xi)}{q_{\xi}(\xi,t)}\right)^{2}.
 \eqno(4.9)$$
By introducing the auxiliary function
$$p(\xi,t) = m(q(\xi,t),t)q_{\xi}(\xi,t)=
\frac{(w^{\prime}_{0}(\xi))^{2}}{q_{\xi}(\xi,t)}\eqno(4.10)$$
and using the formula
$$u(x,t ) = {1\over 2}\int_{\Bbb R} e^{-|x - q(\eta, t)|} m(q(\eta,t),t)
  q_{\eta}(\eta,t)\,d\eta, \eqno(4.11)$$
it follows from (4.4) and (4.9) that
$$\eqalign{
\dot q(\xi,t)& = {1\over 2} \int_{\Bbb R} e^{-|q(\xi,t)- q(\eta, t)|} p(\eta,t)\,d\eta\cr
& = \frac{\delta H}{\delta p}, \cr}\eqno(4.12)$$
and hence
$$\eqalign{
\dot p(\xi,t)& ={1\over 2} p(\xi,t) \int_{\Bbb R} sgn(\xi-\eta)
e^{-|q(\xi,t) - q(\eta, t)|} p(\eta,t)\,d\eta\cr
& = - \frac{\delta H}{\delta q}\cr}\eqno(4.13)$$
where 
$$H = {1\over 4}\int_{\Bbb R^{2}} e^{-|q(\xi,t)- q(\eta, t)|} p(\xi,t)p(\eta,t)\,
  d\xi d\eta.\eqno(4.14)$$
Therefore, (4.12), (4.13) is a Hamiltonian system  with constraint
given by (4.10).(See Remark 4.1 (b) and the appendix for 
more details.) If $q(\xi,t)$, $p(\xi,t)$ satisfy (4.12), (4.13)
and we define
$$K(\xi, \eta;t) = {1\over 2}e^{-{1\over 2}|q(\xi,t)- q(\eta, t)|} 
  \sqrt{p(\xi,t)p(\eta,t)},\eqno(4.15)$$
then a direct calculation shows that $K(\xi, \eta;t)$ evolves
under the integro-differential equation (2.24).  
Note that the kernel $K(\cdot,\cdot\,;t)$ defined in (4.15) is
a positive single-pair kernel.  Moreover, 
$K(\cdot,\cdot\,;t)\in C(\Bbb R^{2})\cap L^{2}(\Bbb R^{2})$
and the corresponding operator is of trace class.(See Remark 4.1 (b).)
\smallskip
\noindent{\bf Remark 4.1.} (a) For the class of initial data which we consider
here, it has been established in \c{R} that the solution $u$ of
the Cauchy problem associated to the CH equation (4.1) satisfies
$u\in C^{\infty}((0,\infty), {\Cal S}(\Bbb R)).$  
\smallskip
\noindent (b) On the other hand, by the ODE for $q(\xi,t)$ in (4.4) (and Remark 4.1 (a))
or otherwise, we know that $q(\cdot,t)\in C^{\infty}(\Bbb R)$ and
hence $p(\cdot,t)$ in also in $C^{\infty}(\Bbb R)$ by (4.10).
Now it follows from (4.12) that
$$\dot q_{\xi}(\xi,t) = \left(-{1\over 2} \int_{\Bbb R} 
sgn(\xi-\eta) e^{-|q(\xi,t) - q(\eta, t)|} p(\eta,t)\,d\eta\right)
q_{\xi}(\xi,t).$$
Solving, we find 
$$q_{\xi}(\xi,t) = \exp\left(-{1\over 2}\int_{0}^{t} \int_{\Bbb R} 
sgn(\xi-\eta) e^{-|q(\xi,\tau) - q(\eta, \tau)|} p(\eta,\tau)\,d\eta\,
d\tau\right). \qquad\,\qquad (**)$$
Since $P = \int_{\Bbb R} p(\eta,t)\,d\eta$ is a conserved quantity \c{C1},
we obtain from (**) that
$e^{-Pt/2}\leq q_{\xi}(\xi, t)\leq e^{Pt/2}$, $\xi\in \Bbb R,\, t>0.$
Consequently, for each $t>0,$ the function $q(\cdot,t)$ is a 
diffeomorphism of the
line which is strictly increasing and satisfies 
$\lim_{\xi\to \pm \infty} q(\xi,t)=\pm \infty$.  
Moreover, $q(\cdot,t)$ is a tempered distribution.
Now, if we continue by
differentiating (**) repeatedly, we can show by
induction that there exist constants $C_{k}(t)>0$ for $t > 0$ such that
$|\partial^{k}_{\xi}q(\xi,t)|\leq C_{k}(t), \,k=2,3,\cdots.$
Since $m(\xi,0)= (w^{\prime}(\xi))^{2}$ is a rapidly decreasing function,
it is now easy to see by using these bounds that
$p(\cdot,t)$ as defined in (4.10) is also a rapidly decreasing function.
\smallskip
\noindent (c) We will show that the map
$(q,p)\mapsto \gamma(q,p)\doteq  
2 e^{-{1\over 2}|q(\xi)- q(\eta)|} \sqrt{p(\xi)p(\eta)}$
is a Poisson map in the appendix. (This is where
the discussion in (b) above will be used in setting
up the domain of the map.)   However, it is
important to point out that the image of this map
is not invariant under $Ad^{*}_{G_R},$ although it
is invariant under the CH-flow. 
\medskip

To solve the Cauchy problem
$$u_t -u_{xxt} + 3uu_x = 2u_{x}u_{xx} + uu_{xxx}, \quad u(x,0) = u_{0}(x)\eqno(4.16)$$
in the space ${\Cal S}(\Bbb R)$ under the above assumptions, we proceed as 
follows.
From the initial data $u_{0}$ which satisfies assumptions (i) and (ii)
above, we obtain the initial conditions for $q(\xi,t)$ and
$p(\xi,t)$:
$$q(\xi,0) = \xi, \quad p(\xi,0) = m(\xi,0) = (w^{\prime}_{0}(\xi))^{2}. 
\eqno(4.17)$$
From this, we obtain the kernel $K(\xi,\eta)$ of the operator $\bk_{0}$,
namely,
$$K(\xi,\eta) = {1\over 2} e^{-{1\over 2}|\xi -\eta|}w^{\prime}_{0}(\xi)
  w^{\prime}_{0}(\eta).\eqno(4.18)$$
\smallskip
\noindent{\bf Remark 4.2.} Note that from our assumption on the initial 
data $u_0$ and (4.18) above, it is clear that $\bk_{0}\in \frak p_{*}.$  Hence 
it follows from Proposition 2.8 that the solution $\bk(t)$ of
(3.22) is on the coadjoint orbit
${\Cal O}_{\bk_{0}} =\{\,Ad^{*}_{G_R}(g)\bk_{0}\mid g\in G\,\}\subset \frak p_{*}$
for all $t.$   
\medskip
Now we solve the initial value problem (3.22) whose solution
$\bk(t)$ is given by the formula in (3.21) and we denote
the kernel corresponding to $\bk(t)$ by $K(\xi,\eta,\,;t).$
Since $K(\xi, \eta,\,;t)$ is given by the formula in (4.15)
where $q(\xi,t)$, $p(\xi,t)$ are solutions of (4.12), (4.13)
with initial condition given in (4.17), we can recover
$p(\eta, t)$ from the formula
$$p(\eta, t) = 2 K(\eta,\eta;t). \eqno(4.19)$$
From (4.15) and (4.19), it follows that
$$ e^{-|q(\xi,t) -q(\eta,t)|} = \frac{K(\xi,\eta;t)^2}{K(\xi,\xi;t)K(\eta,\eta;t)}.
  \eqno(4.20)$$
Hence we obtain
$$\dot q(\eta, t) = \int_{\Bbb R} \frac{K(\eta,\zeta;t)^2}{K(\eta,\eta;t)}\,
  d\zeta \eqno(4.21)$$
and consequently, we can determine $q(\eta,t)$ from the formula
$$q(\eta,t) = \eta + \int_{0}^{t} \left(\int_{\Bbb R} \frac{K(\eta,\zeta;\tau)^2}
 {K(\eta,\eta;\tau)}\,d\zeta \right) d\tau.\eqno(4.22)$$
Finally, we obtain the solution of the Cauchy problem (4.16) from
$$u(x,t) = {1\over 2}\int_{\Bbb R} e^{-|x - q(\eta, t)|} p(\eta,t)\,d\eta
 \eqno(4.23)$$
where $p(\eta,t)$ and $q(\eta,t)$ are given in (4.19) and (4.22)
above.

Note that the explicit formula for $K(\xi,\eta;t)$ is given 
in (3.36) where in the present case we can interpret equality
in the pointwise sense as the kernels of all operators involved
in (3.35) are continuous.  Alternatively, we can make use of
Mercer's theorem in writing down the explicit formula for
$K(\xi,\eta;t)$.  In this connection, note that $\bk_{0}$
is compact and hence its spectrum $\sigma(\bk_{0})$ is discrete
with no limit points except possibly at $0.$  Moreover,
$\bk_{0}$ has no negative eigenvalues.  To see this,
let $\lambda$
be an eigenvalue of $\bk_{0}$ and $\phi$ a corresponding normalized
eigenfunction.  Then from
$${1\over 2} \int_{\Bbb R}  e^{-{1\over 2}|\xi -\eta|}w^{\prime}_{0}(\xi)
  w^{\prime}_{0}(\eta) \phi(\eta)\, d\eta = \lambda \phi(\xi),\eqno(4.24)$$
we have
$$\lambda = {1\over 2} \int_{\Bbb R} \left(\int_{\Bbb R}
e^{-{1\over 2}|\xi-\eta|}(w^{\prime}_{0}\phi)(\eta)\,d\eta \right)
(w^{\prime}_{0}\phi)(\xi)\,d\xi.\eqno(4.25)$$
Therefore, on applying the Parseval formula and the convolution theorem
in Fourier transforms to the right hand of (4.25), we obtain\footnotemark\
\footnotetext{We owe this argument and the sharpening of (4.28)
to the referee.}%
$$\lambda = {1\over 4\pi} \int_{\Bbb R} \frac{|(\widehat{w^{\prime}_{0}\phi})(k)|^{2}}{k^{2} + {1\over 4}}\, dk \eqno(4.26)$$
from which the positivity is clear. (Here $\widehat{w^{\prime}_{0}\phi}$ is 
the Fourier transform of $w^{\prime}_{0}\phi.$)
Hence
this shows that $\bk_{0}$ has no negative eigenvalues.  To give
an estimate of  $\sigma(\bk_{0})$, note that the integral on
the right hand side of (4.25) obeys the inequalities 
$$\eqalign{&
\left|\int_{\Bbb R} \left( \int_{\Bbb R} e^{-{1\over 2}|\xi -\eta|}
w^{\prime}_{0}(\eta) \phi(\eta)\, d\eta \right) 
w^{\prime}_{0}(\xi) \phi(\xi)
\, d\xi \right|\cr
\leq & \left(\int_{\Bbb R} w^{\prime}_{0}(\xi)|\phi(\xi)|\,d\xi\right)^2\cr
\leq & \int_{\Bbb R} (w^{\prime}_{0}(\xi))^2 \,d\xi =  P,\cr}\eqno(4.27)$$
where we have used the Cauchy-Schwarz inequality and the
normalization of $\phi$ in going from
the second line to the last one in (4.27) above.
Combining the above analysis,  we can now conclude that
$$\sigma(\bk_{0}) \subset [0, P/2]. \eqno(4.28)$$
Let $\sigma(\bk_{0}) = \{\lambda_{i}\}_{i=1}^{\infty}$ and let 
$\{\phi_{i}\}_{i=1}^{\infty}$ be the corresponding normalized
eigenfunctions.  Then we have
$$\bk_{0} = \sum_{i=1}^{\infty}  \lambda_{i} \phi_{i}\otimes \phi_{i}\eqno(4.29)$$
so that from (3.21), (3.20) and the orthogonality of $b_{+}(t)$, we find
$$\eqalign{
\bk(t) & = \sum_{i=1}^{\infty} \lambda_{i} (b_{+}(t)^{-1}\phi_{i})\otimes
          (b_{+}(t)^{-1}\phi_{i})\cr
       & = \sum_{i=1}^{\infty} \lambda_{i} e^{-t\lambda_{i}}(b_{-}(t)^{-1}\phi_{i})\otimes
          (b_{-}(t)^{-1}\phi_{i}).\cr}\eqno(4.30)$$
Since $\sigma(\bk(t))= \sigma(\bk_{0})$, 
the series
$\sum_{i=1}^{\infty} \lambda_{i} e^{-t\lambda_{i}}(b_{-}(t)^{-1}\phi_{i})(\xi)
          (b_{-}(t)^{-1}\phi_{i})(\eta)$
converges absolutely and uniformly on compact sets and it follows
from (4.30) and Mercer's theorem that
$$K(\xi,\eta;t) = \sum_{i=1}^{\infty} \lambda_{i} e^{-t\lambda_{i}}(b_{-}(t)^{-1}\phi_{i})
(\xi)(b_{-}(t)^{-1}\phi_{i})(\eta)\eqno(4.31)$$
where from (3.31) and (3.32),
$$(b_{-}(t)^{-1}\phi_{i})(x) = \phi_{i}(x) - \int_{-\infty}^{x}
  \left(S(y,x;t) +\frac {(D_{2}(\bs(t)|_{(-\infty, x)})S(\cdot,x;t))(y)}
    {{\det}_{2}\left(e^{-t\bk_{0}}|_{(-\infty, x)}\right)}\right)\phi_{i}(y)\,dy.
   \eqno(4.32)$$
\smallskip
\noindent{\bf Remark 4.3.} Under different assumption on $m(\xi,0)$
and using an entirely different method, the Lagrangian form of 
the CH equation was integrated in terms of certain Fredholm
determinants in \c{M}.  In particular, the spectral problem
in \c{M} is given by 
$$({1\over 4}-\partial^2_{x})f(x,t) = \lambda m(x,t) f(x,t).\eqno(4.33)$$
On the other hand, our spectral problem reads:
$${1\over 2}\int_{-\infty}^{\infty}  e^{-{1\over 2}|q(\xi,t)- q(\eta, t)|} 
  \sqrt{p(\xi,t)p(\eta,t)} \phi(\eta,t)\,d\eta = \lambda \phi(\xi,t).
\eqno(4.34)$$
It is a natural question to ask if (4.34) can be derived from
(4.33).   To put it differently, can we derive  the
kernel (4.15) (as discovered in \c{C1}) from (4.33)?
As the referee pointed out to us, the answer is an
affirmative yes.  Indeed, under the assumption in \c{M},
(4.33) can be rewritten as 
$$\int_{-\infty}^{\infty} e^{-{1\over 2}|x-y|} m(y,t)f(y,t)\,dy = {f(x,t)\over \lambda}.
\eqno(4.35)$$
Then the change to the Lagrangian variable $x=q(\xi,t)$ gives
$$\int_{-\infty}^{\infty} e^{-{1\over 2}|q(\xi,t)-q(\eta,t)|} p(\eta,t)f(q(\eta,t),t)\,d\eta
 = {f(q(\xi,t),t)\over \lambda}\eqno(4.36)$$
where we have invoked the definition of $p(\cdot,t)$ in (4.10).
For our class of initial data, $p(\xi,t)>0$ for all $t.$
Therefore,if we multiply both sides of (4.36) by $\sqrt{p(\xi,t)}$,
and setting $\phi(\xi,t) = f(q(\xi,t),t)\sqrt{p(\xi,t)},$ the
result is (4.34) above (modulo a factor $1\over 2$) provided we
replace $\lambda$ by ${1\over \lambda}$.

\bigskip
\bigskip

\newpage

\head
{\bf Appendix}
\endhead
\bigskip

Let $U = {\Cal S}^{\prime}(\Bbb R)\cap \{\hbox{strictly increasing 
diffeomorphisms}\,\, q:\Bbb R\longrightarrow \Bbb R \}$, where
${\Cal S}^{\prime}(\Bbb R)$ is the space of tempered distributions,
and let $V$ be the subset of ${\Cal S}(\Bbb R)$ consisting of those
$p\in {\Cal S}(\Bbb R)$ which are strictly positive.
We equip $U\times V$ with the Poisson bracket
$$\{{\Cal F}_1, {\Cal F}_2\}(q,p)= \int_{\Bbb R} \left(
{\frac{\partial {\Cal F}_1}{\partial q(\xi)}}
{\frac{\partial {\Cal F}_2}{\partial p(\xi)}} -
{\frac{\partial {\Cal F}_1}{\partial p(\xi)}}
{\frac{\partial {\Cal F}_2}{\partial q(\xi)}}\right)\,d\xi.\eqno(A1)$$
For $(q,p)\in U\times V$, put
$$\lp(\xi) = \sqrt{2 p(\xi)} e^{{1\over 2} q(\xi)},\,\,\, \lm(\xi)=
\sqrt{2 p(\xi)} e^{-{1\over 2} q(\xi)}.\eqno(A2)$$

\proclaim
{Proposition A}  The map $\gamma:U\times V\longrightarrow \fg^*_{R}\simeq
\fg_{R}$ defined by
$$\gamma(q,p) \doteq 2 \,e^{-{1\over 2}|q(\xi)- q(\eta)|} \sqrt{p(\xi)p(\eta)}=
                \cases \lp(\xi)\lm(\eta), 
                & \xi\leq \eta \\
                 \lp(\eta)\lm(\xi), & \xi>\eta, \endcases\eqno(A3)$$
is a Poisson map.
\endproclaim

\demo
{Proof} We want to show
$$\{F_{1}\circ \gamma, F_{2}\circ \gamma \} (q,p) =
  \{F_1,F_2\}_{R}(\gamma (q,p)).$$
To compute the Fr\'echet derivatives, we proceed as follows.
First of all,
$$\int_{\Bbb R} {\frac{\partial (F_{i}\circ\gamma)}{\partial q(\xi)}}
  \widetilde q(\xi)\,d\xi =\left(dF_{i}(\gamma(q,p)), {d\over d\epsilon}
  {\Big|_{\epsilon =0}}
  \gamma (q + \epsilon \widetilde q, p)\right).\qquad (\dagger)$$
By using the definition of $\gamma,$ a straight forward computation
shows that
$$
{d\over d\epsilon}{\Big|_{\epsilon =0}}\gamma (q + \epsilon \widetilde q, p)
\doteq \cases {1\over 2}\lp(\xi)\lm(\eta)(\widetilde q(\xi)-
                 \widetilde q(\eta)),& \xi\leq \eta,\\
                {1\over 2}\lp(\eta)\lm(\xi)(\widetilde q(\eta)-
                 \widetilde q(\xi)),& \xi>\eta.\endcases$$
Substitute this into ($\dagger$), it follows after some
manipulation that
$${\frac{\partial (F_{i}\circ\gamma)}{\partial q(\xi)}} =
{1\over 2} \lp(\xi)((\Pi_{\frak l}\, dF_{i}(\g))^{*}\lm)(\xi)
 -\lm(\xi)((\Pi_{\frak l}\, dF_{i}(\g))\lp)(\xi)$$
where $dF_{i}(\g)$ is the shorthand for $dF_{i}(\g(q,p))$ which 
we will use from now onwards.  In a similar way, 
we can show that
$$
{d\over d\epsilon}{\Big|_{\epsilon =0}}\gamma (q, p + \epsilon \widetilde p)
\doteq \cases {\lm(\eta)\over \lm(\xi)}\widetilde p(\xi) + 
               {\lp(\xi)\over \lp(\eta)}
               \widetilde p(\eta), & \xi\leq\eta,\\
               {\lm(\xi)\over \lm(\eta)}\widetilde p(\eta) + 
               {\lp(\eta)\over \lp(\xi)}
               \widetilde p(\xi), & \xi>\eta ,\endcases$$
and
$${\frac{\partial (F_{i}\circ\gamma)}{\partial p(\xi)}} =
{1\over \lp(\xi)}((\Pi_{\frak l}\, dF_{i}(\g))\lp)(\xi) +
{1\over \lm(\xi)}((\Pi_{\frak l}\, dF_{i}(\g))^{*}\lm)(\xi).$$
Next, we substitute the Fr\'echet derivatives into the
expression for the Poisson bracket between $F_{1}\circ \gamma$
and $F_{2}\circ \gamma,$
this gives
$$\aligned
&\{F_{1}\circ \gamma, F_{2}\circ \gamma \} (q,p)\\
= & {1\over 2} \int_{\Bbb R} [\,((\Pi_{\frak l}\, dF_{1}(\g))^{*} \lm)(\xi)
    ((\Pi_{\frak l}\, dF_{2}(\g))\lp)(\xi)\\ 
  & \quad \,\,\,- ((\Pi_{\frak l}\, dF_{1}(\g)) \lp)(\xi)
    ((\Pi_{\frak l}\, dF_{2}(\g))^{*} \lm)(\xi)\\ 
  & \quad \,\,\, - (F_{1}\leftrightarrow F_{2})\,]\,d\xi\\ 
= & \int_{\Bbb R} [\,((\Pi_{\frak l}\, dF_{1}(\g))^{*}\lm)(\xi)
    ((\Pi_{\frak l}\, dF_{2}(\g))\lp)(\xi)\\ 
  & \quad \,\,\,- ((\Pi_{\frak l}\, dF_{1}(\g)) \lp)(\xi)
    ((\Pi_{\frak l}\, dF_{2}(\g))^{*} \lm)(\xi)\,]\,d\xi.
\endaligned
$$
But now by using the fact that $\Pi_{\frak l}\, dF_{1}(\g),$
 $\Pi_{\frak l}\, dF_{2}(\g)\in \frak l,$
we can show that
$$\int_{\Bbb R} \,((\Pi_{\frak l}\, dF_{1}(\g))^{*}\lm)(\xi)
    ((\Pi_{\frak l}\, dF_{2}(\g))\lp)(\xi)\, d\xi =
   (\g(q,p), \Pi_{\frak l}\,dF_{1}(\g)\circ \Pi_{\frak l}\,dF_{2}(\g)).$$
Interchanging the indices $1$ and $2$, we obtain
$$\int_{\Bbb R} \,((\Pi_{\frak l}\, dF_{2}(\g)^{*} \lm)(\xi)
    ((\Pi_{\frak l}\, dF_{1}(\g)\lp)(\xi)\, d\xi =
   (\gamma(q,p), \Pi_{\frak l}\,dF_{2}(\g)\circ \Pi_{\frak l}\,
   dF_{1}(\g)).$$
Therefore, when we subtract the second expression from the
first, the result is
$$\{F_{1}\circ \gamma, F_{2}\circ \gamma \} (q,p) =
   (\g(q,p), [\Pi_{\frak l}\, dF_{1}(\g), \Pi_{\frak l}\, dF_{2}(\g)]). \qquad 
(\ddagger)$$
To complete the proof, note that 
$(\g(q,p), [\Pi_{\frak k}\, dF_{1}(\g), \Pi_{\frak k}\, dF_{2}(\g)]) =0$
as $\gamma(q,p)\in \frak p$.
Consequently, when we add 
$-(\g(q,p), [\Pi_{\frak k}\, dF_{1}(\g), \Pi_{\frak k}\, dF_{2}(\g)])$
to the right hand side of  $(\ddagger)$, the assertion follows.
\pf
\enddemo

\enddocument